\documentclass[11pt,a4paper,reqno]{amsproc}
\usepackage{longtable,cite}
\usepackage{lscape}
\usepackage{multirow}
\usepackage[pdftex,bookmarks=true]{hyperref}

\theoremstyle{definition}

\theoremstyle{remark}

\numberwithin{equation}{section}
\usepackage{makecell}
\usepackage{geometry} 
\usepackage{caption}
\usepackage{pdflscape}
\usepackage{stmaryrd}
\usepackage{rotating}
\usepackage{subcaption}  
\usepackage{multirow} 
\usepackage{graphicx} 
\usepackage{array} 
\usepackage{algorithm}
\usepackage{amsmath}
\usepackage{tcolorbox}
\usepackage{caption}
\usepackage{ragged2e}
\usepackage{lipsum}  
\usepackage{placeins}
\usepackage{amsmath, amssymb}
\usepackage{tikz}
\usepackage{makecell}
\usetikzlibrary{shapes.geometric, arrows}
\tikzstyle{startstop} = [rectangle, rounded corners, minimum width=3cm, minimum height=1cm,text centered, draw=black, fill=white]
\tikzstyle{process} = [rectangle, minimum width=3cm, minimum height=1cm, text centered, draw=black, fill=white]
\tikzstyle{decision} = [diamond, minimum width=3cm, minimum height=1cm, text centered, draw=black, fill=white]
\tikzstyle{arrow} = [thick,->,>=stealth]

\usepackage{float}
\usepackage{blindtext}
\usepackage{graphicx}
\graphicspath{{Diagrams/}}
\usepackage{csquotes}
\usepackage{adjustbox}
\usepackage{soul}
\usepackage{xcolor} 
\usepackage{hyperref} 

\usepackage{float}
\usepackage{fancyhdr}
\title [{Assessing Finite Element Choice in Structural Topology Optimization and A Posteriori Error Estimation}]{Assessing Finite Element Choice in Structural Topology Optimization and A Posteriori Error Estimation} 
\author{Jyotiranjan Nayak}
\address{Department of Mathematics, SRM University AP, Andhra Pradesh 522240, India}
\email{jyotiranjan\_n@srmap.edu.in}
\author{Shafeequdheen P}
\address{Department of Mathematics, SRM University AP, Andhra Pradesh 522240, India}
\email{shafeequdheen\_p@srmap.edu.in}

\author{Vijayakrishna rowthu}
\address{Department of Mathematics, SRM University AP, Andhra Pradesh 522240, India}
\email{vijayakrishna.r@srmap.edu.in}

\date{}

\usepackage{cite}

\begin{document}

\maketitle	
	
	\begin{abstract}
			
This study investigates the impact of finite element selection on structural topology optimization using the SIMP (Solid Isotropic Material with Penalization) method. Specifically, it compares linear (P1) and quadratic (P2) triangular elements with the conventional bi-linear quadrilateral (Q1) elements. Numerical experiments performed on benchmark problems including a cantilever beam, a bridge structure, and a beveled beam reveal notable differences in both the final optimized objective value (compliance) and the accuracy of the finite element solutions. The accuracy is evaluated using an a posteriori error estimator, highlighting the influence of element type on solution quality and optimization performance.

	\end{abstract}

\section{Introduction}

In recent decades, the topological optimization has established itself as a fundamental and versatile methodology for the design of efficient and innovative structures across a broad range of engineering disciplines. The origins of topology optimization trace back to A.G.M. Michell’s pioneering work in 1904\cite{Michell1904}, where he established theoretical criteria for minimum-weight truss structures. After a long dormancy, renewed interest came in the 1972s\cite{Rozvany1972a, Rozvany1972b} when G. Rozvany and others extended Michell’s concepts to continuum structures, leading to the formulation of the “Optimal Layout Theory” by Prager \& Rozvany \cite{PragerRozvany1977}. A major breakthrough occurred in 1988 when Bendsøe and Kikuchi\cite{Bendsoe1988} introduced a homogenization-based numerical method, enabling optimal material distribution in continuum domains and marking the beginning of modern topology optimization. Since then, several key methodologies have emerged, notably the Solid Isotropic Material with Penalization (SIMP) \cite{Bendsoee1989, bendsoe2001topology, mlejnek1992some}, Evolutionary Structural Optimization (ESO) \cite{xie1993simple, xie1997basic, xie2010recent} and its extension Bi-directional Evolutionary Structural Optimization (BESO) \cite{zhuang2023172, querin1998evolutionary, querin2000computational}, and the level set method \cite{allaire2002level, wang2003level, allaire2004structural}, each contributing unique perspectives and advantages to the field. 
 
  These methodologies have enabled topology optimization to become a central tool in high-performance structural design. Optimized layouts are increasingly applied in additive manufacturing, aerospace and automotive design, architectural structures, and biomedical engineering. For instance, topology-optimized parts are now routinely used in lightweight UAV frames, lattice-structured orthopedic implants, heat exchangers, and artistic structures requiring both functional and aesthetic performance. Recent trends involve integration with machine learning, uncertainty quantification, and real-time design systems, marking a new frontier in intelligent, manufacturable, and robust design under diverse conditions.

Numerical methods for topology optimization have evolved significantly since the seminal work of Bendsøe and Kikuchi \cite{Bendsoe1988}, which laid the groundwork for a range of density-based approaches. Among these, SIMP method, first introduced by Bendsøe \cite{Bendsoee1989} and later refined by Zhou and Rozvany \cite{Zhou1991}, has become one of the most widely used techniques. SIMP discretizes the design domain into finite elements, assigning each element a material density variable, typically ranging from 0 (void) to 1 (solid). To promote discrete solutions, a penalization scheme is applied to intermediate densities, typically using a penalty exponent \( p \geq 3 \), which encourages convergence toward a black-and-white design with few or no intermediate densities. This formulation is effective for a variety of applications, including multi-material design, eigenvalue optimization, thermo-elastic and fluid-structure problems, and additive manufacturing constraints. A key advancement that contributed to SIMP’s widespread use was the development of Sigmund's 99-line MATLAB code \cite{Sigmund2001}, which provided an accessible platform for educational and prototyping purposes, making SIMP an indispensable tool for researchers and engineers alike. In the SIMP framework, the design domain is fixed, and the material distribution is optimized by updating the element-wise densities. However, while SIMP offers a continuous approach to topology optimization, it often results in semi-dense, intermediate elements that may not be manufacturable. To address this, various regularization techniques have been proposed, such as filtering (via classical sensitivity or density filters, projection techniques, morphology-based filters, or Helmholtz-type filters), as well as geometric constraint techniques, including perimeter and gradient constraints \cite{Bendsoee1989, Sigmund2001}. These methods alleviate issues such as gray areas, checkerboard patterns, and mesh dependency, which arise from the ill-posed nature of unconstrained topology optimization problems. The final solution in SIMP-based topology optimization is typically obtained through iterative optimization algorithms, such as the Optimality Criteria (OC) method \cite{Bendsoe1989}, the Methods of Moving Asymptotes (MMA) algorithm \cite{Sigmund2001}, or other mathematical programming-based techniques. These methods, in conjunction with regularization schemes, ensure more stable and manufacturable solutions. Through these advancements, SIMP has maintained its position as the leading approach for topology optimization across a wide range of engineering disciplines. 

In addition to the widely adopted SIMP method, the ESO technique originally introduced by Xie and Steven in the 1990s~\cite{xie1993simple, xie1996evolutionary, xie1997basic}, although similar ideas had been explored earlier has emerged as a popular alternative for topology optimization, especially in industrial applications. ESO is a heuristic, hard-kill approach that begins with a fully solid design domain and iteratively removes inefficient material based on local rejection criteria. This binary element-based formulation ($\chi \in \{0, 1\}$) produces clear, black-and-white topologies but often suffers from convergence challenges and sensitivity to initial configurations, which can lead to local optima. To overcome these limitations, the method was extended into the BESO framework~\cite{yang1999bidirectional, querin2000computational, sun2011topology}, which permits both removal and addition of material. This enhancement improves the method’s capacity to explore the design space more effectively and helps avoid poor local minima. New material is typically introduced near regions with high sensitivity or within voids, estimated via interpolation of the displacement field. However, such extrapolations may be inconsistent with the modeling of solid regions, and hard-kill versions of BESO can exhibit poor convergence behavior. To address this, Huang and Xie~\cite{xie2010recent} introduced modifications incorporating stabilization techniques, mesh-independent filters, and the use of historical design data, significantly improving robustness and convergence over time. More recently, advanced BESO methods have incorporated body-fitted meshes and regularization techniques to enhance boundary smoothness and numerical stability. For instance, Zhuang et al.~\cite{zhuang2023172} developed a body-fitted BESO approach based on nonlinear diffusion regularization, which leads to sharper boundaries and improved convergence. Furthermore, they introduced TriTOP172, a compact 172-line MATLAB implementation, providing an efficient and educational tool for topology optimization on body-fitted meshes. Alternatively, soft-kill approaches retain void regions as low-density elements rather than eliminating them, enabling the evaluation of sensitivities across the entire design domain. Such methods, including those proposed by Zhu et al.~\cite{zhu2007bi} and Huang and Xie~\cite{huang2009bi}, combine the BESO framework with SIMP-style penalized density interpolation, thereby improving numerical stability and extending the applicability of ESO-based methods.

Another important class of methods is based on the level-set framework. Originally developed by Osher and Sethian \cite{osher1988fronts} for tracking moving interfaces, the level-set method was later adapted for topology optimization by Allaire et al. \cite{allaire2004structural}, Wang et al. \cite{wang2003level}, and others. In this approach, the structural boundary is represented implicitly as the zero level set of a scalar function. This representation allows smooth boundary evolution and naturally accommodates topological changes such as merging components or nucleation of holes. Compared with density-based formulations, level-set methods generate sharp structural boundaries and avoid intermediate material densities.

Various formulations of level-set based topology optimization have been developed, including approaches based on shape derivatives, topological derivatives, and phase-field models \cite{wang2003level, allaire2004structural}. However, these methods often require carefully chosen initial configurations containing small holes throughout the design domain to enable sufficient topological evolution. To overcome this limitation, hole nucleation techniques (sometimes referred to as bubble methods) have been introduced to allow the creation of new void regions during the optimization process. Additional improvements such as Hamilton–Jacobi evolution equations, reinitialization procedures, and regularization strategies are commonly employed to ensure numerical stability and convergence \cite{allaire2004structural}.

While significant progress has been made in the development of topology optimization algorithms, the accuracy of the underlying finite element discretization remains a crucial factor influencing the quality and reliability of the optimized structures. In this context, a posteriori error estimation provides an effective framework for assessing the accuracy of computed finite element solutions. Calculations always require reliable control over the accuracy of approximations obtained, and the development of analytical and practical tools for such error control constitutes the main purpose of a posteriori error estimation analysis. Various approaches to derive estimates for elliptic-type boundary value problems, with errors measured in the energy norm, have been suggested by many authors \cite{ainsworth1997posteriori, babuvvska1978error, zienkiewicz1987simple}. Residual-based estimators, which measure element residuals and flux jumps across element interfaces, provide reliable indicators of discretization error and play an important role in adaptive finite element methods. A comprehensive overview of these techniques is given by Verfürth \cite{verfurth1999review}. Within topology optimization, such estimators offer valuable insight into how discretization errors influence the performance and reliability of optimized structural layouts.

In topology optimization, the choice of finite elements significantly influences both the optimized structural layout and the accuracy of the finite element solution. Conventional SIMP-based implementations commonly employ bilinear quadrilateral elements ($Q_1$), which remain computationally efficient but may exhibit higher discretization errors and reduced stress-field accuracy, particularly in problems involving complex load transfer. In contrast, triangular discretizations using linear ($P_1$) and quadratic ($P_2$) elements provide greater geometric flexibility, enabling improved representation of evolving structural features during optimization. The numerical results obtained for the cantilever beam, bridge structure, and beveled beam demonstrate that triangular elements consistently yield lower residual-based \textit{a posteriori} error estimates than $Q_1$ elements, indicating superior solution accuracy. Among the triangular formulations, $P_1$ elements produce the lowest compliance values in most cases, suggesting stiffer optimized designs while maintaining moderate computational cost. Quadratic $P_2$ elements deliver smoother displacement and stress fields due to higher-order interpolation, but this improvement comes with increased computational expense and does not necessarily lead to lower compliance. Overall, the study confirms that triangular discretizations outperform traditional quadrilateral elements in terms of accuracy and optimization performance across different problem configurations. The results further indicate that $P_1$ elements provide an effective balance between solution quality and computational efficiency, making them particularly suitable for large-scale topology optimization problems, whereas $P_2$ elements are advantageous when higher solution smoothness is required.

The remainder of this paper is organized as follows. Section 2 presents the necessary preliminaries, including the strong and weak formulations of the governing equations together with the fundamentals of residual-based a posteriori error estimation. Section 3 describes the mathematical formulation of the topology optimization problem and the adopted optimization strategy. Section 4 reports the numerical experiments and provides a comparative analysis of the results obtained with different finite element discretizations. Finally, Section 5 summarizes the main conclusions of the study.

\section{Preliminaries}	
	The strong form of the  linear elasticity problem seeks a displacement field $\mathbf{u} : \Omega \to \mathbb{R}^d$ satisfying the following boundary value problem (see \cite{ern2004theory}):
\begin{align}
	\left\{
	\begin{aligned}
		-\nabla \cdot \sigma (\mathbf{u}) &= \mathbf{f}, && \text{in } \Omega, \\
		\mathbf{u} &= \mathbf{0}, && \text{on } \Gamma_D, \\
		\sigma(\mathbf{u}) \cdot \hat{n} &= \mathbf{g}, && \text{on } \Gamma_N,
	\end{aligned}
	\right.
	\label{1.1}
\end{align}
where $\mathbf{f}$ is the body force per unit volume, $\mathbf{g}$ is the prescribed surface traction, and $\hat{n}$ is the outward unit normal to the boundary. Further the constitutive equations are given by

\begin{align}
	\left\{
	\begin{aligned}
		\sigma(\mathbf{u}) &= A \varepsilon(\mathbf{u}) = \lambda\, \text{tr}(\boldsymbol{\varepsilon}(\mathbf{u}))\, \mathbf{I} + 2\mu\, \boldsymbol{\varepsilon}(\mathbf{u}), \\
		\boldsymbol{\varepsilon}(\mathbf{u}) &= \frac{1}{2} \left( \nabla \mathbf{u} + \nabla \mathbf{u}^\top \right).
	\end{aligned}
	\right.
	\label{1.2}
\end{align}
The stress tensor $\boldsymbol{\sigma}$ and the linearized strain tensor $\boldsymbol{\varepsilon}$ are related through Hooke’s law for isotropic materials (see \cite{johnson2009numerical}),
where $\lambda$ and $\mu$ are the Lamé parameters, $\text{tr}(\cdot)$ denotes the trace operator, and $\mathbf{I}$ is the identity tensor. The elasticity tensor $A$ is used in more general anisotropic formulations (see \cite{ciarlet1997mathematical}).

\subsection{Weak formulation}
Assuming the solution $\textbf{u}$ is made up of a linear combination of some test functions from $V = \left\{ \textbf{v} \in [H^1(\Omega)]^d : \textbf{v}|_{\Gamma_D} = \mathbf{0} \right\}$. We refer to Lemma~2.1.11 and Lemma~2.1.12 in the lecture notes by Emma Cinatl  ~\cite{cinatl2018finite} for the necessary integration-by-parts identities used in the derivation of the weak formulation.

The weak form of the linear elasticity problem can be compactly written as the variational equation
\begin{equation}
\int_\Omega \sigma(\textbf{u}): \varepsilon(\textbf{v}) = \int_{\Omega} \mathbf{f} \cdot \textbf{v} dx + \int_{\Gamma_N} \mathbf{g} \cdot \textbf{v} ds.
	\label{eq:main}
\end{equation}
\begin{equation}
	a(\textbf{u}, \textbf{v}) = L(\textbf{v})
	\label{eq:auvl}
\end{equation}
where the bilinear form $a(\cdot,\cdot)$ and the linear form $L(\cdot)$ are defined as in equation~\eqref{eq:a} and ~\eqref{eq:l} respectively.
\begin{equation}
	a(\textbf{u}, \textbf{v})   = \int_\Omega \sigma(\textbf{u}): \varepsilon(\textbf{v}) = \int_{\Omega} \lambda (\nabla \cdot \textbf{u})(\nabla \cdot \textbf{v}) dx 
	+ \int_{\Omega} 2\mu \boldsymbol{\varepsilon}(\textbf{u}) : \boldsymbol{\varepsilon} (\textbf{v})  dx
	\label{eq:a}
\end{equation}

\begin{equation}
	L(\textbf{v}) = \int_{\Omega} \mathbf{f} \cdot \textbf{v} dx + \int_{\Gamma_N} \mathbf{g} \cdot \textbf{v} ds.
	\label{eq:l}
\end{equation}

Setting \(\textbf{u} = \sum_{j=1}^{N} C_j \phi_j\) and \(\textbf{v} = \phi_i\), where $\{\phi_i\}_{i=1}^{N}$
are the nodal basis functions that span the finite-dimensional space \(V_h \subset V\) and equation~\eqref{eq:auvl} can be written as system of linear equations by using equation~\eqref{1.2} as follows
\begin{equation}
	\sum_{j=1}^{N} \left[ \lambda \int_{\Omega} (\nabla \cdot {\phi}_i)(\nabla \cdot {\phi}_j) \, dx 
	+ 2\mu \int_{\Omega} {\boldsymbol{\varepsilon}}({\phi}_i) : {\boldsymbol{\varepsilon}}({\phi}_j) \, dx \right] C_j 
	= \int_{\Omega} \mathbf{f} \cdot {\phi}_i \, dx 
	+ \int_{\Gamma_N} \mathbf{g} \cdot {\phi}_i \, ds.
	\label{eq:weak-9}
\end{equation}

Further, the equation~\eqref{eq:weak-9} can be represented in the matrix form as  \( \mathbf{K}\mathbf{U} = \mathbf{F} \), where, \( K_{ij} = a(\phi_i, \phi_j)\), \(F_i = L(\phi_i) \) and \(\mathbf{U}\) consists of the unknowns \( C_j \), \( 1 \leq i, j \leq n \). The well-posedness and proof details can be found in Emma Cinatl \cite{cinatl2018finite}.

The finite element discretization of the problem ~\eqref{eq:auvl} is given by 
\begin{equation}
\int_\Omega \sigma(u_h): \varepsilon(v_h) = 	\int_{\Omega} f \cdot v_h dx + \int_{\Gamma_N} g \cdot v_h ds.
	\label{eq:dis}
\end{equation}	
	
\subsection{A Posteriori Error Analysis}
A posteriori error estimation provides a practical framework for evaluating the accuracy of finite element solutions based on the computed approximation itself. By measuring residuals of the governing equations within elements and jumps across element interfaces, these estimators deliver reliable indicators of discretization error and are widely used in adaptive finite element methods \cite{verfurth1999review}.

Let $\Omega \subset \mathbb{R}^d$ ($d=2,3$) be a bounded domain with boundary 
$\partial \Omega = \Gamma_D \cup \Gamma_N$, where $\Gamma_D \cap \Gamma_N = \emptyset$. 
Consider a mesh $\mathcal{T}_h$ partitioning $\Omega$ into elements $K$, each with diameter $ h_K := \mathrm{diam}(K)$, 
and let $h_e := \|e\|$ denote the length of an edge $e$.

Denote the set of all edges by $ \mathcal{E}_h = \mathcal{E}_h^{\mathrm{int}} \cup \mathcal{E}_h^{N} \cup \mathcal{E}_h^{D}$, where $\mathcal{E}_h^{\mathrm{int}}$ is the set of interior edges, and 
$\mathcal{E}_h^{N}$ and $\mathcal{E}_h^{D}$ denote edges lying on the Neumann and Dirichlet boundaries, respectively.

For an interior edge $e \in \mathcal{E}_h^{\mathrm{int}}$, shared by elements $K^+$ and $K^-$, 
the jump of a sufficiently smooth vector field $\mathbf{w}$ across $e$ is defined as

\begin{equation}
	[\![ \mathbf{w} ]\!] := \mathbf{w}^{+} \cdot \mathbf{n}^{+} + \mathbf{w}^{-} \cdot \mathbf{n}^{-},
	\label{eq:jump_definition}
\end{equation}

where $\mathbf{n}^{\pm}$ are the outward unit normals to elements $K^{\pm}$ on the edge $e$. 
On boundary edges, $\mathbf{n}$ denotes the outward unit normal to $\Omega$.

Let $\mathbf{u} \in V$ and $\mathbf{u}_h \in V_h$ be the exact and finite element solutions, respectively, 
where $V$ is the appropriate Sobolev space and $V_h \subset V$ is its finite-dimensional subspace. 
Define the error

\begin{equation}
	\boldsymbol{\zeta} := \mathbf{u} - \mathbf{u}_h .
	\label{eq:error_definition}
\end{equation}

The residual functional $R : V \rightarrow \mathbb{R}$ is defined by

\begin{equation}
	R(v) := L(v) - a(\mathbf{u}_h , v) = a(\boldsymbol{\zeta}, v),
	\quad \forall v \in V ,
	\label{eq:residual_functional}
\end{equation}

where $a(\cdot,\cdot)$ is the bilinear form corresponding to the weak formulation of the elasticity problem, and $L(\cdot)$ represents the linear functional associated with the external loads.

Applying integration by parts elementwise, the residual can be expressed as

\begin{equation}
	R(v) =
	\sum_{K \in \mathcal{T}_h} \int_K R_K v \, dx
	+
	\sum_{e \in \mathcal{E}_h^{\mathrm{int}}} \int_e J_e v \, ds
	+
	\sum_{e \in \mathcal{E}_h^{N}} \int_e R_e^{N} v \, ds ,
	\label{eq:residual_decomposition}
\end{equation}

where the local residual terms are defined as

\begin{equation}
	R_K := \mathbf{f} + \nabla \cdot \sigma(\mathbf{u}_h), 
	\qquad
	J_e := [\![ \sigma(\mathbf{u}_h)\cdot \mathbf{n}]\!], 
	\qquad
	R_e^{N} := \mathbf{g} - \sigma(\mathbf{u}_h)\cdot \mathbf{n}.
	\label{eq:local_residuals}
\end{equation}

Here, $\mathbf{f}$ denotes the body force, $\mathbf{g}$ the prescribed Neumann boundary data, 
and $\sigma(\mathbf{u}_h)$ is the stress tensor evaluated from the finite element solution.

Choosing $v = \boldsymbol{\zeta}$ yields the energy norm identity

\begin{equation}
	\| \boldsymbol{\zeta} \|_a^2 := a(\boldsymbol{\zeta},\boldsymbol{\zeta}) = R(\boldsymbol{\zeta}).
	\label{eq:energy_identity}
\end{equation}

Standard a posteriori error estimation theory provides the reliability bound

\begin{equation}
	\| \boldsymbol{\zeta} \|_a^2 
	\lesssim 
	\sum_{K \in \mathcal{T}_h} h_K^2 \| R_K \|_{0,K}^2
	+
	\sum_{e \in \mathcal{E}_h^{\mathrm{int}}} h_e \| J_e \|_{0,e}^2
	+
	\sum_{e \in \mathcal{E}_h^{N}} h_e \| R_e^{N} \|_{0,e}^2 ,
	\label{eq:reliability_bound}
\end{equation}

where $\|\cdot\|_{0,D}$ denotes the $L^2$-norm over the domain $D$.

Define the elementwise error indicator

\begin{equation}
	\eta_K^2 :=
	h_K^2 \| R_K \|_{0,K}^2
	+
	\sum_{e \subset \partial K \cap \mathcal{E}_h^{\mathrm{int}}}
	h_e \| J_e \|_{0,e}^2
	+
	\sum_{e \subset \partial K \cap \Gamma_N}
	h_e \| R_e^{N} \|_{0,e}^2 ,
	\label{eq:local_indicator}
\end{equation}

and the global error estimator

\begin{equation}
	\eta :=
	\left(
	\sum_{K \in \mathcal{T}_h} \eta_K^2
	\right)^{1/2}.
	\label{eq:global_estimator}
\end{equation}

	\section{SIMP For Topology Optimization Problem}
	 A widely used method for this purpose is the SIMP approach, originally proposed by Bendsøe \cite{Bendsoee1989} and later formalized with Sigmund \cite{Bendsoee2001}. In SIMP, each finite element is assigned a design variable $x_e \in [0, 1]$, representing its relative material density.
	The SIMP for topology optimization problem is then formulated as the minimization of compliance (i.e., maximization of stiffness) under a volume constraint:
	\begin{align}
		\left\{
		\begin{aligned}
			\min_{\mathbf{x}} \quad & c(\mathbf{x}) = \mathbf{U}^\top \mathbf{K}(\mathbf{x}) \mathbf{U} = \sum_{e=1}^{N} x_e^p\, \mathbf{u}_e^\top \mathbf{K}_0 \mathbf{u}_e \\
			\text{subject to} \quad & \dfrac{V(\mathbf{x})}{V_0} = f, \quad \mathbf{K}(\mathbf{x}) \mathbf{U} = \mathbf{F},\quad 0 \leq x_{\min} \leq x_e \leq 1 \quad \text{for } e = 1, \ldots, N.
		\end{aligned}
		\right.
		\label{eq:simp}
	\end{align}
	\noindent
	Here, \( \mathbf{x} = \{x_1, x_2, \ldots, x_{N}\}^\top \) is the vector of design variables. The term \( V(\mathbf{x}) = \displaystyle \sum_{e=1}^{N} x_e v_e \) denotes the total material volume, where \( v_e \) is the volume of element \( e \), and \( V_0 \) is the volume of the full design domain. The volume fraction \( f \in (0,1) \) defines the allowed proportion of material. The small lower bound \( x_{\min} \) ensures numerical stability and prevents singularities in the stiffness matrix \( \mathbf{K} \). The equilibrium behavior of an elastic structure, derived from the  weak form of the governing equations (see the appendix for the derivation), is discretized using finite element methods and expressed as $ \mathbf{K} \, \mathbf{U} = \mathbf{F}$. This forms the foundation for structural analysis in topology optimization, where the objective is to find an optimal material distribution within a prescribed domain satisfying the volume constraint. The stiffness of element \( e \) is interpolated using a penalized power-law model
	\begin{equation}
		\mathbf{K}_e(x_e) = x_e^p \mathbf{K}_0,
	\end{equation}
	where \( \mathbf{K}_0 \) is the stiffness matrix of a fully solid element and \( p > 1 \) is the penalization exponent, typically set to 3, which promotes near-binary designs by penalizing intermediate densities. This interpolation makes the global stiffness matrix \( \mathbf{K}(\mathbf{x}) \) explicitly dependent on the design vector \( \mathbf{x} \), thereby coupling the design and analysis steps.

	\subsection{Finite Element Analysis (FEA)}
	
	In the SIMP-based topology optimization framework, the finite element method (FEM) is used to evaluate structural responses based on the equilibrium equation \( \mathbf{K} \mathbf{U} = \mathbf{F} \). While prior studies commonly used bi-linear quadrilateral (Q1) elements, this work employs linear (\(P_1\)) and quadratic (\(P_2\)) triangular elements for both design space and solution space discretizations. The associated shape functions serve as basis functions for interpolating displacement fields within each element.
	
	\begin{figure}[htbp]
		\centering
		\begin{subfigure}[htbp]{0.36\textwidth}
			\centering
			\textbf{Shape functions for \(P_1\) elements}
			
			\[
			\begin{aligned}
				\phi_1(x, y) &= 1 - x - y, \\
				\phi_2(x, y) &= x, \\
				\phi_3(x, y) &= y.
			\end{aligned}
			\]
			
			\begin{tikzpicture}[scale=1.6]
				\draw[thick] (0,0) -- (1,0) -- (0,1) -- cycle;
				\filldraw[black] (0,0) circle (1pt);
				\filldraw[black] (1,0) circle (1pt);
				\filldraw[black] (0,1) circle (1pt);
				\node[below left] at (0,0) {$\phi_1$};
				\node[below right] at (1,0) {$\phi_2$};
				\node[above right] at (0,1) {$\phi_3$};
			\end{tikzpicture}
			
			\caption{Master element with \(P_1\) basis functions.}
			\label{fig:P1-shape-functions}
		\end{subfigure}
		\hfill
		\begin{subfigure}[htbp]{0.63\textwidth}
			\centering
			\textbf{Shape functions for \(P_2\) elements}
			
			\[
			\begin{aligned}
				\phi_1 &= (1 - x - y)(1 - 2x - 2y), \\
				\phi_2 &= x(2x -1), \quad
				\phi_3 = y(2y -1), \\
				\phi_4 &= 4x(1 - x - y), \\
				\phi_5 &= 4xy, \quad
				\phi_6 = 4y(1 - x - y).
			\end{aligned}
			\]
			
			\begin{tikzpicture}[scale=1.5]
				\draw[thick] (0,0) -- (1,0) -- (0,1) -- cycle;
				\filldraw[black] (0,0) circle (1pt);
				\filldraw[black] (1,0) circle (1pt);
				\filldraw[black] (0,1) circle (1pt);
				\filldraw[black] (0.5,0) circle (1pt);
				\filldraw[black] (0.5,0.5) circle (1pt);
				\filldraw[black] (0,0.5) circle (1pt);
				\node[below left] at (0,0) {$\phi_1$};
				\node[below right] at (1,0) {$\phi_2$};
				\node[above right] at (0,1) {$\phi_3$};
				\node[below] at (0.5,0) {$\phi_4$};
				\node[right] at (0.5,0.5) {$\phi_5$};
				\node[left] at (0,0.5) {$\phi_6$};
			\end{tikzpicture}
			
			\caption{Master element with \(P_2\) basis functions.}
			\label{fig:P2-shape-functions}
		\end{subfigure}
	\end{figure}

	The strain-displacement matrix \( [B] \), constructed from the spatial derivatives of the shape functions, enables computation of strain as \( \boldsymbol{\varepsilon} = [B] \, \mathbf{u}_e \), and stress via the constitutive relation \( \boldsymbol{\sigma} = A \boldsymbol{\varepsilon} \), where \( A \) is the elasticity matrix. For linear elements, \( [B] \) is constant over each element, whereas for quadratic elements it varies spatially within each element. The elemental strain energy density (SED) is expressed as
	\[
	\text{SED}_e = \int_e \mathbf{u}_e^\top [B]^\top A [B] \mathbf{u}_e dV
	\]
	and the corresponding stiffness matrix is \( \mathbf{K}_e = \int_{V^e} [B]^\top A [B] \, dV \). For linear elements, this is efficiently approximated by scaling with the element area. To ensure numerical accuracy, especially with quadratic elements, a 7-point Gaussian quadrature rule is used for numerical integration over each triangular element
	\[
	\int_{A_e} f(x, y) \, dx\,dy \approx A_e \sum_{i=1}^{7} w_i \, f(x_i, y_i).
	\]
	
\begin{figure}[h]
		\centering
		\begin{minipage}[htbp]{0.48\textwidth} 
			\begin{tcolorbox}[colback=gray!5,colframe=black!50,
				title=\textbf{Pseudocode for the loop}, width=\textwidth,
				left=2mm, right=2mm, top=1mm, bottom=1mm,
				boxsep=1mm, halign=left]
				\hspace*{1em}~Initialize $x$ with volfrac for each element;\\ \hspace*{1em}~ $\texttt{loop}=0$, $\texttt{rchange}=1$\\
				~\textbf{while} $\texttt{rchange}>0.01$ \textbf{do}\\
				\hspace*{1em}~$\texttt{loop}\leftarrow\texttt{loop}+1$;  $x_{\text{old}}\leftarrow x$\\
				\hspace*{1em}~Compute displacement $\mathbf{U}\leftarrow\text{FEAnalysis}$\\
				\hspace*{1em}~$c\leftarrow0$\\
				\hspace*{1em}~\textbf{for} each element $e$ \textbf{do}\\
				\hspace*{2em}Compute $\mathbf{K}_e$, $\mathbf{U}_e$, 
				$\text{comp}\!\leftarrow\!\mathbf{U}_e^T\mathbf{K}_e\mathbf{U}_e$\\
				\hspace*{2em}$c\!\leftarrow\!c+(x_e)^{\texttt{penal}}\text{comp}$; \\ 
				\hspace*{2em}$dc_e\!\leftarrow\!-\texttt{penal}(x_e)^{\texttt{penal}-1}\text{comp}$\\
				\hspace*{1em}~\textbf{end for}\\
				\hspace*{1em}~$dc\leftarrow\text{SensitivityFilter}$; \\
				\hspace*{1em}~$x\leftarrow\text{OptimalityCriteriaUpdate}$\\
				\hspace*{1em}~$\texttt{rchange}\leftarrow
				\dfrac{\max(|x-x_{\text{old}}|)}{\max(x_{\text{old}})}$;  \\
				\hspace*{1em}~$\text{Plot the Updated Design}$\\
				~\textbf{end while}
			\end{tcolorbox}
		\end{minipage}
		\hfill
		\begin{minipage}[htbp]{0.48\textwidth} 
			\centering
			\resizebox{\textwidth}{!}{ 
				\begin{tikzpicture}[node distance=1.5cm, scale=0.35, every node/.style={scale=1}]
					\node (start) [startstop] {\textbf{Start}};
					\node (domain) [process, below of=start, align=center] {%
						\textbf{Domain Initialization}\\
						(Nodes, Elements, Material Properties)};
					\node (BC) [process, below of=domain, align=center]{\textbf{FE Analysis}\\ (Apply BCs, Calc $K_e$, Solve $KU = F$)};
					\node (objective) [process, below of= BC] {\textbf{Compute Objective} $C(x_e)$};
					\node (sensitivity) [process, below of=objective] {\textbf{Compute Sensitivities} $\dfrac{\partial C}{\partial x_e}$};
					\node (filtering) [process, below of=sensitivity] {\textbf{Apply Filter}};
					\node (update) [process, below of=filtering] {\textbf{Update Design}};
					\node (plot) [process, right of=update, xshift = 5cm] {\textbf{Plot Updated Design}};
					\node (loop) [decision, below of=update, yshift=-1cm] {\textbf{Converged?}};
					\node (end) [startstop, below of=loop, yshift=-1.2cm] {\textbf{End}};
					\draw [arrow] (start) -- (domain);
					\draw [arrow] (domain) -- (BC);
					\draw [arrow] (BC) -- (objective);
					\draw [arrow] (objective) -- (sensitivity);
					\draw [arrow] (sensitivity) -- (filtering);
					\draw [arrow] (filtering) -- (update);
					\draw [arrow] (update) -- (plot);
					\draw [arrow] (update) -- (loop);
					\draw [arrow] (loop.south) -- ++(0,-0.1) node[midway, right] {Yes} -- (end.north);
					\draw [arrow] (loop.west) -- ++(-11,0) coordinate(tmp1)
					-- ++(0,24.3) coordinate(tmp2)
					-- (BC.west) node[midway, above left] {No};
				\end{tikzpicture}
			}
		\end{minipage}
		\caption{Left: Pseudocode; Right: Flowchart of the SIMP method}
		\label{fig:pseudocode_flowchart}
	\end{figure}
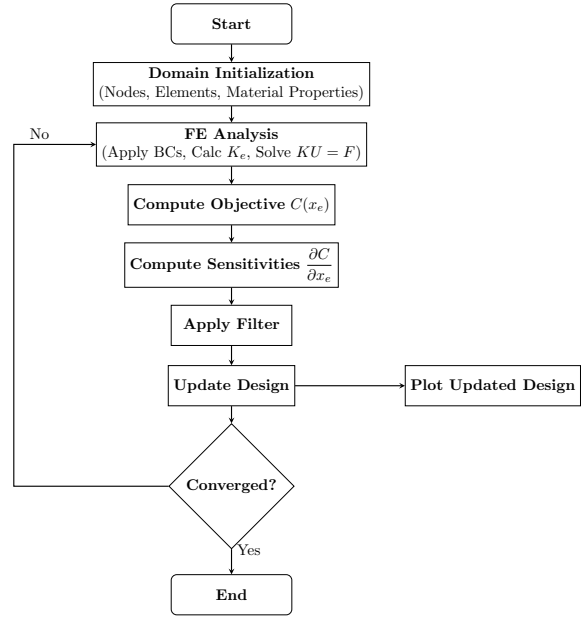

	Finally, the global stiffness matrix \( \mathbf{K} \) is assembled by summing contributions from all elements using standard sparse matrix techniques
	\[
	\mathbf{K} = \sum_{e=1}^{N_e} \mathbf{P}_e^\top \mathbf{K}_e \mathbf{P}_e,
	\]
	where \( \mathbf{P}_e \) maps local to global degrees of freedom. This FEA infrastructure directly supports the SIMP-based optimization loop by providing displacement fields and elemental energies for updating the design variables.
 The schematic representation of the methods that have been used in this model to solve the optimization problem has been provided in Figure~\ref{fig:pseudocode_flowchart} with the pseudocode for the loop.
	
Building upon the well-established SIMP framework coupled with finite element analysis, the present study continues to serve as a foundation for investigating the influence of finite element choice on both solution quality and optimization outcomes. By integrating linear ($P_1$) and quadratic ($P_2$) triangular elements alongside the conventional bilinear quadrilateral ($Q_1$) elements into the optimization process, we systematically quantify differences in shapes viacompliance and discretization error via a posteriori error estimation. This approach reveals the critical impact of element selection not only on the optimized layout but also on the accuracy and reliability of the finite element solutions. The novelty of these findings lies in highlighting how finite element type directly affects optimization performance, thereby providing valuable insights for the development of new algorithmic strategies. 

\section{Numerical Results and Discussion}

In this section, numerical experiments are conducted on three benchmark problems, namely the cantilever beam, bridge structure, and beveled beam, in order to evaluate the influence of element type on both the optimized compliance and the accuracy of the finite element solution. The present work can be viewed as an extension of the classical SIMP-based topology optimization framework popularized by Ole Sigmund \cite{Sigmund2001}, where bilinear quadrilateral ($Q_1$) elements are typically employed for finite element discretization. In the well-known 99-line topology optimization code and related studies, $Q_1$ elements were used to demonstrate the effectiveness and simplicity of the SIMP approach. In the present study, the same SIMP formulation is retained, but the finite element discretization is extended to include triangular elements, namely linear ($P_1$) and quadratic ($P_2$) elements. This extension enables a comparative investigation of how different element types influence the optimized structural layouts, compliance values, and the accuracy of the finite element solution when assessed using a residual-based a posteriori error estimator. In all cases, identical material properties, boundary conditions, loading conditions, volume fraction, and penalization parameters are used to ensure a consistent comparison. The setup for each of these test cases, including detailed material properties, boundary conditions, and optimization parameters, is provided below.

For the computation of the a posteriori error estimator, the volume fraction is taken as $volfrac = 1$, corresponding to a fully solid design domain. When smaller volume fractions (e.g., $40\%$ or $50\%$ of the initial design domain) are used during topology optimization, the resulting structures contain thin members and solid--void interfaces that introduce strong stress gradients and large residual contributions in the error estimator. As a result, the estimated error can become significantly large during the early stages of the optimization process and gradually decreases as the topology evolves toward the final optimized configuration. To obtain a stable and meaningful comparison of discretization errors across different element types, the error estimation is therefore performed for the full material case.

\vspace{0.5cm}			
			
\begin{figure}[htbp]
		\centering
		\renewcommand{\arraystretch}{1.5}
		\begin{tabular}{|
				>{\centering\arraybackslash}m{5cm}|
				>{\centering\arraybackslash}m{5cm}|
				>{\centering\arraybackslash}m{5cm}|}
			\hline 
			\textbf{$Q_1$} & \textbf{$P_1$}& \textbf{$P_2$} \\

			\hline
			
			\includegraphics[width=5cm,height=4cm]{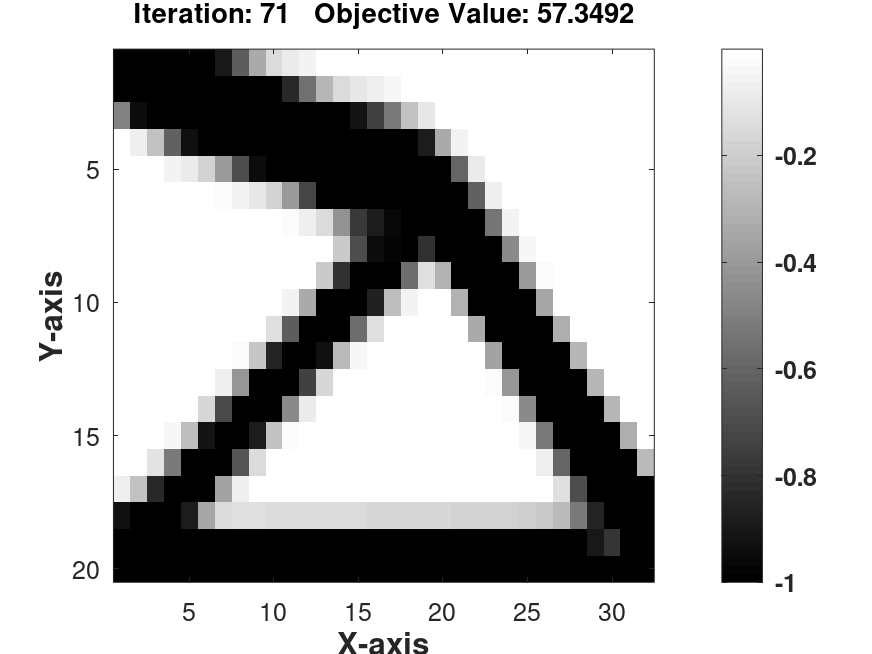} 
				Elements = 640  & 
			 \includegraphics[width=5cm,height=4cm]{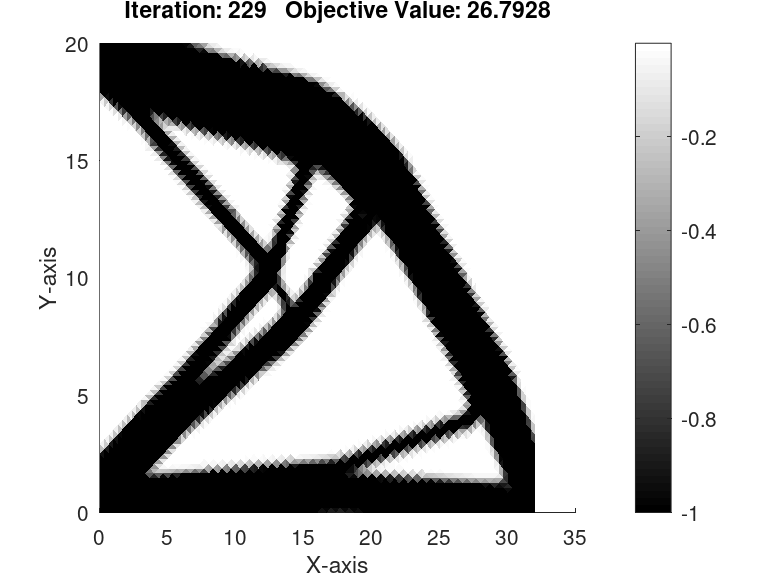} 
			 	Elements = 9216 &  
             \includegraphics[width=5cm,height=4cm]{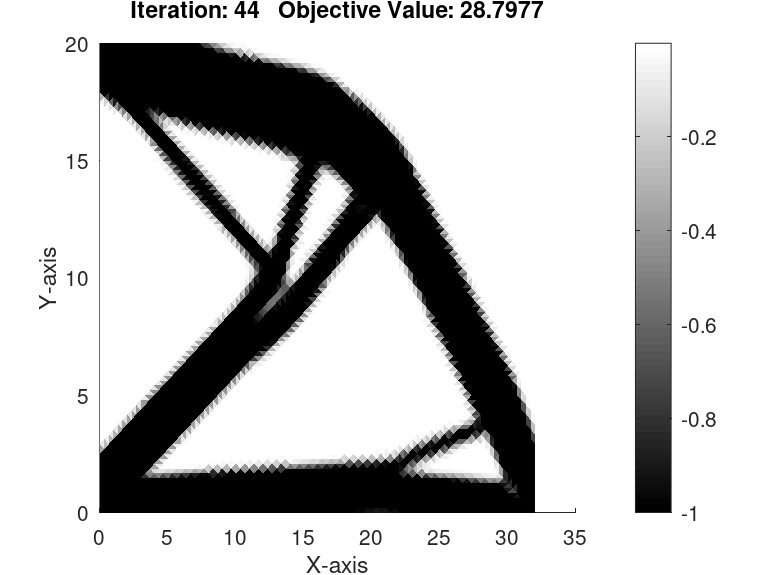}  
             	Elements = 9216  \\
			\hline
			& 
			 \includegraphics[width=5cm,height=4cm]{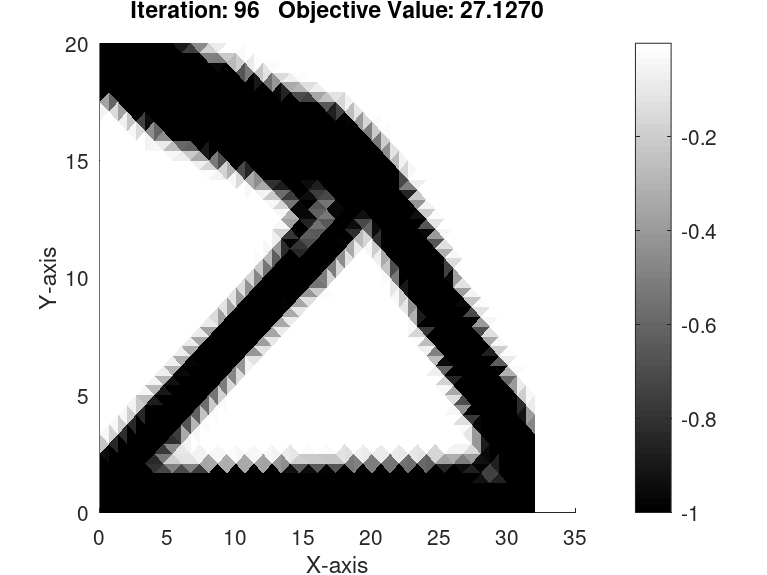} 
			 	Elements = 2304 &
			 \includegraphics[width=5cm,height=4cm]{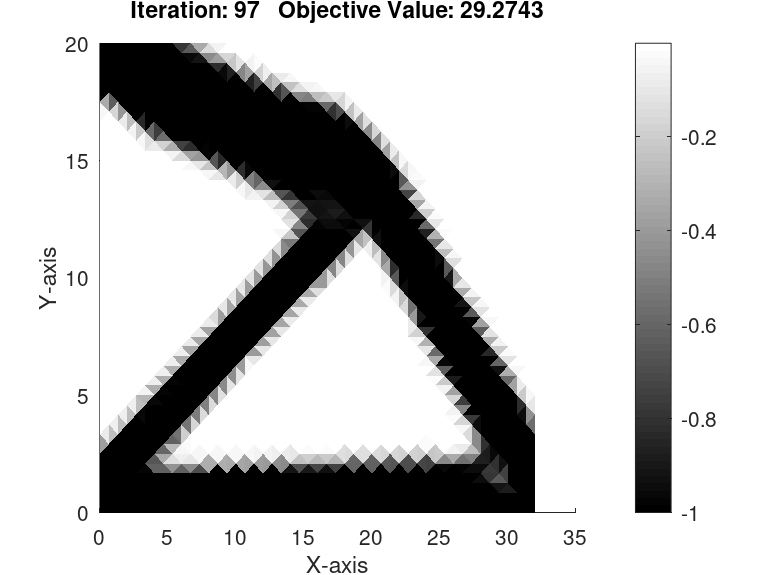}
			 	Elements = 2304  \\
			 \hline			
		\end{tabular}
		\caption{Comparison for \textbf{Cantilever Beam  of Dimension $32 \times 20$}}
		\label{tab:cc1}
	\end{figure}

		\begin{minipage}[t]{0.38\textwidth} 
			\centering
			\begin{tikzpicture}[scale=1.0]
				\fill[gray!50] (0,0) rectangle (4,2);
				\draw[thick] (0,0) rectangle (4,2);
				\foreach \y in {0.1,0.3,...,1.9}
				\draw[black] (-0.2,\y) -- (0,\y+0.1);
				\filldraw[red] (4,0) circle (0.05);
				\draw[very thick, ->, red] (4.2,0) -- (4.2,-0.4);
				\node[below] at (4.3,-0.4) {\footnotesize Force};
				\node[rotate=90] at (-0.5,1) {\footnotesize Fixed Nodes};
			\end{tikzpicture}
			\captionof{figure}{\textbf{Cantilever Beam} with fixed and force node. $E = 1$, $\nu = 0.3$, $volfrac = 0.4$, $penal = 3$, and $\vec{F} = (0, -1)$.}
			\label{fig:cbeam}
		\end{minipage}
		\hfill
	\begin{minipage}[h]{0.5\textwidth} 
		\textbf{Cantilever Beam:} Figure~\ref{fig:cbeam} illustrates the details of the cantilever beam, where the fixed nodes are located entirely along the west side of the domain. A force is applied at the bottom corner on the east side, directed along the negative Y-axis. The parameters used are described within the same figure. Figure~\ref{tab:cc1} \& ~\ref{tab:cc2} presents the optimized topologies obtained for different resolutions of the grid. 
	\end{minipage}

	\begin{table}[htbp]
		\centering
		\renewcommand{\arraystretch}{1.2}
		
		\begin{tabular}{|
				>{\centering\arraybackslash}m{3cm}|
				>{\centering\arraybackslash}m{3cm}|
				>{\centering\arraybackslash}m{3cm}|
				>{\centering\arraybackslash}m{3cm}|}
			\hline 
			
			\textbf{Element Type} & 
			\textbf{No. of Elements} & 
			\textbf{Final Objective Value} & 
			\textbf{Iterations} \\
			\hline
			
			Q1 & 640 & 57.3492 & 71 \\
			\hline
			
			\multirow{2}{*}{P1} & 9216 & 26.7928 & 229 \\
			\cline{2-4}
			& 2304 & 27.1270 & 96 \\
			\hline
			\multirow{2}{*}{P2} & 9216 & 28.7977 & 44 \\
			\cline{2-4}
			& 2304 & 29.2743 & 97 \\
			\hline
		\end{tabular}
		\caption{Cantilever Beam of Domain dimension $32 \times 20$ }
		\label{tab:1}
	\end{table}

	\begin{figure}[htbp]
	\centering
	\renewcommand{\arraystretch}{1.5}
	\begin{tabular}{|
			>{\centering\arraybackslash}m{5cm}|
			>{\centering\arraybackslash}m{5cm}|
			>{\centering\arraybackslash}m{5cm}|}
		\hline 
\textbf{$Q_1$} & \textbf{$P_1$}& \textbf{$P_2$} \\
		\hline
	\includegraphics[width=5cm,height=4cm]{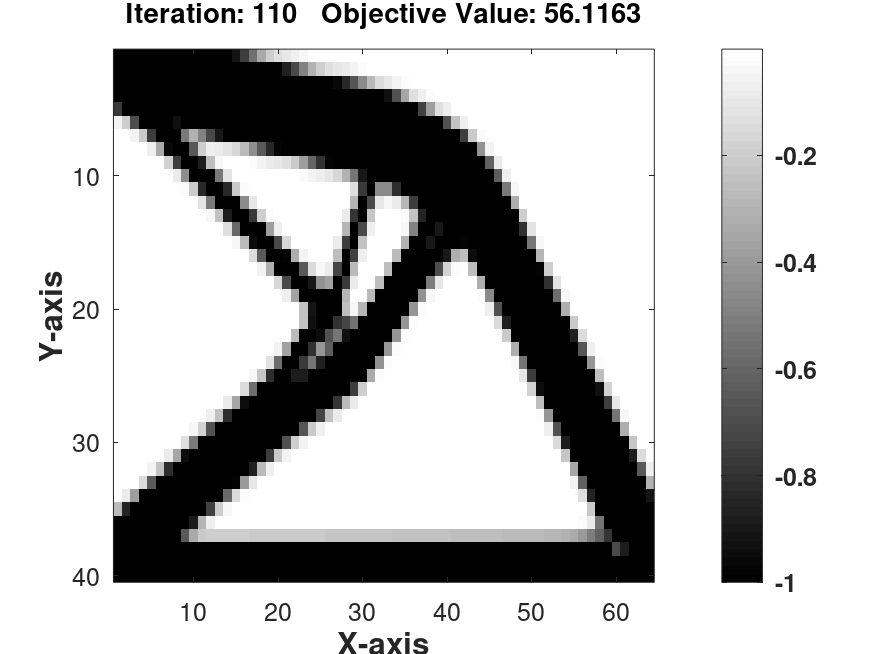}
	Elements = 2560 & 
	\includegraphics[width=5cm,height=4cm]{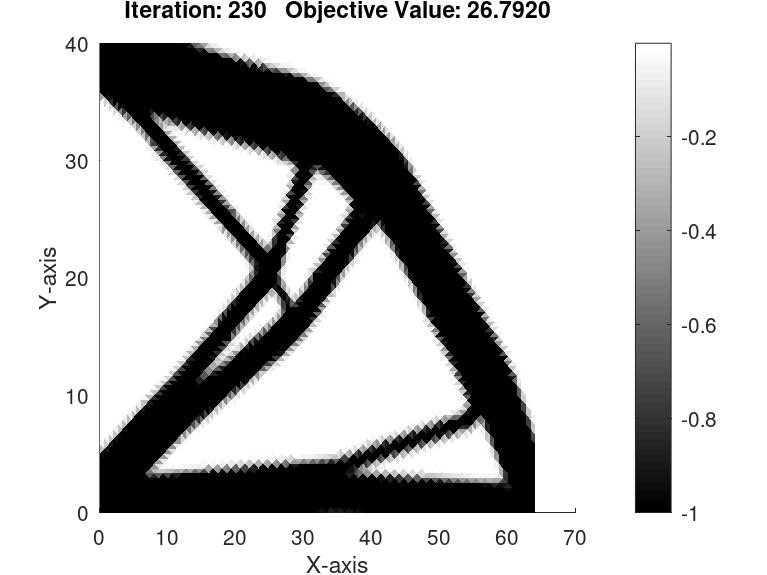}
	Elements = 9216 & 
	\includegraphics[width=5cm,height=4cm]{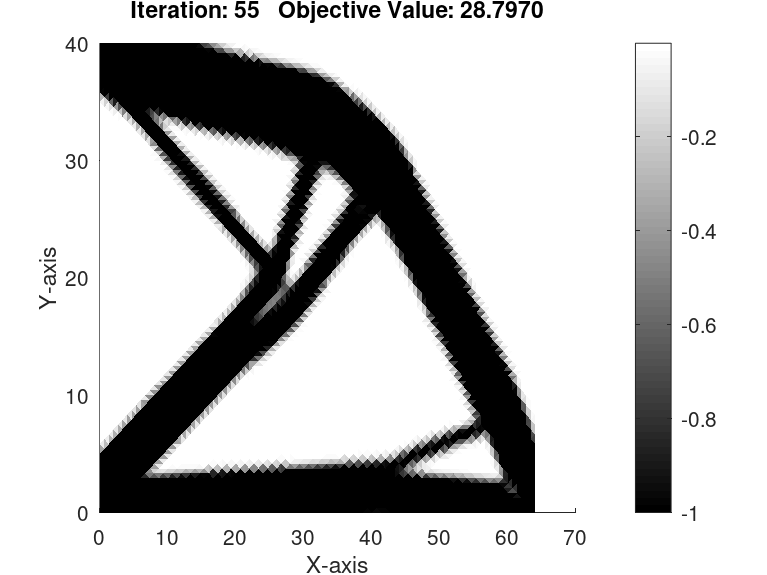}
	Elements = 9216 \\
	\hline
	
	& 
	\includegraphics[width=5cm,height=4cm]{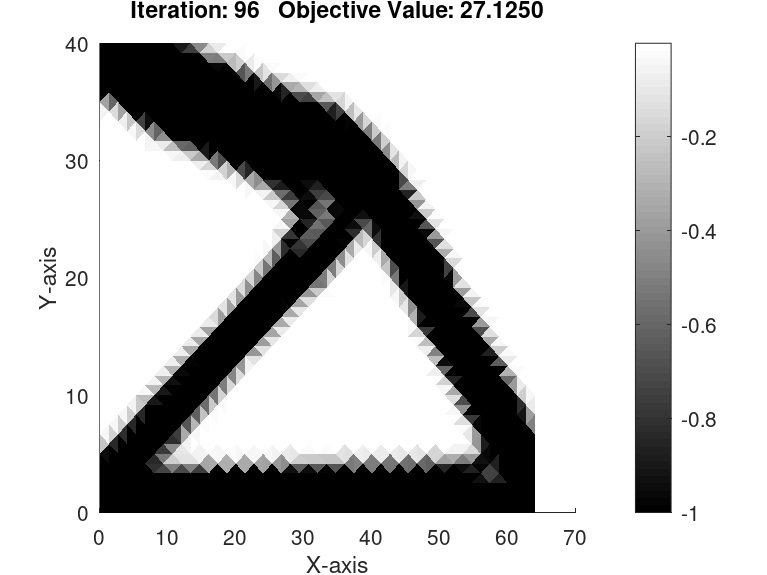} 
	Elements = 2304 & 
	\includegraphics[width=5cm,height=4cm]{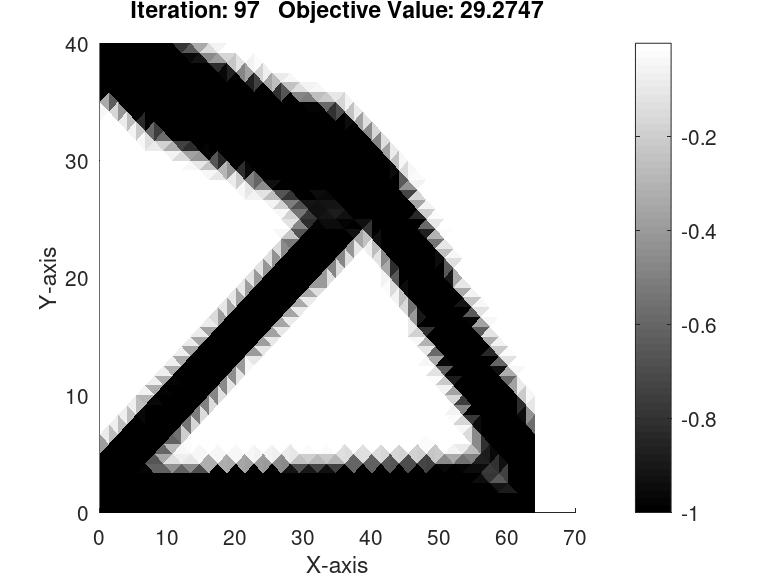}
		Elements = 2304 \\
	\hline

	\end{tabular}
	\caption{Comparison for \textbf{Cantilever Beam of Domain Dimension $64 \times 40$}}
	\label{tab:cc2}
\end{figure}

\begin{table}[htbp]
	\renewcommand{\arraystretch}{1.2}
	
	\begin{tabular}{|
			>{\centering\arraybackslash}m{1.4cm}|
			>{\centering\arraybackslash}m{1.5cm}|
			>{\centering\arraybackslash}m{1.4cm}|
			>{\centering\arraybackslash}m{1.7cm}|
			>{\centering\arraybackslash}m{1.5cm}|
			>{\centering\arraybackslash}m{1.5cm}|
			>{\centering\arraybackslash}m{1.6cm}|
			>{\centering\arraybackslash}m{1.2cm}|
			>{\centering\arraybackslash}m{1.2cm}|}
		\hline
		
		\textbf{Element Type} &
		\textbf{No. of Elements} &
		\textbf{Final Objective Value} &
		\textbf{Iterations} &
		\textbf{Bulk Residual} &
		\textbf{Internal Jump Residual} &
		\textbf{Neumann Residual} &
		\textbf{Local Error $\eta^2$} &
		\textbf{Global Error $\eta$} \\
		
		\hline
		
		Q1 & 2560 & 56.1163 & 110 & 8.1648 & 20.914 & 0.4563 & 29.535 & 5.4346 \\
		\hline
		
		\multirow{2}{*}{P1} 
		& 9216 & 26.7920 & 230 & 0 & 1.2770 & 1.1814 & 2.4585 & 1.5679 \\
		\cline{2-9}
		& 2304 & 27.1250 & 96 & 0 & 1.3909 & 2.3087 & 3.6996 & 1.9234 \\
		\hline
		
		\multirow{2}{*}{P2} 
		& 9216 & 28.7970 & 55 & 4.4235 & 1.7528 & 0.1291 & 6.3053 & 2.5110 \\
		\cline{2-9}
		& 2304 & 29.2747 & 97 & 18.042 & 1.7998 & 0.4462 & 20.288 & 4.5013 \\
		\hline
		
	\end{tabular}
	
	\caption{Comparison of objective values, optimization iterations, and residual-based a posteriori error estimates for the cantilever beam problem on a $64 \times 40$ domain using different finite element discretizations.}
	\label{tab:cantilever_results}
	
\end{table}

\begin{figure}[htbp]
	\centering
	\renewcommand{\arraystretch}{1.5}
	\begin{tabular}{|
			>{\centering\arraybackslash}m{5cm}|
			>{\centering\arraybackslash}m{5cm}|
			>{\centering\arraybackslash}m{5cm}|}
		\hline 
\textbf{$Q_1$} & \textbf{$P_1$}& \textbf{$P_2$} \\
\hline
		
		\includegraphics[width=5cm,height=4cm]{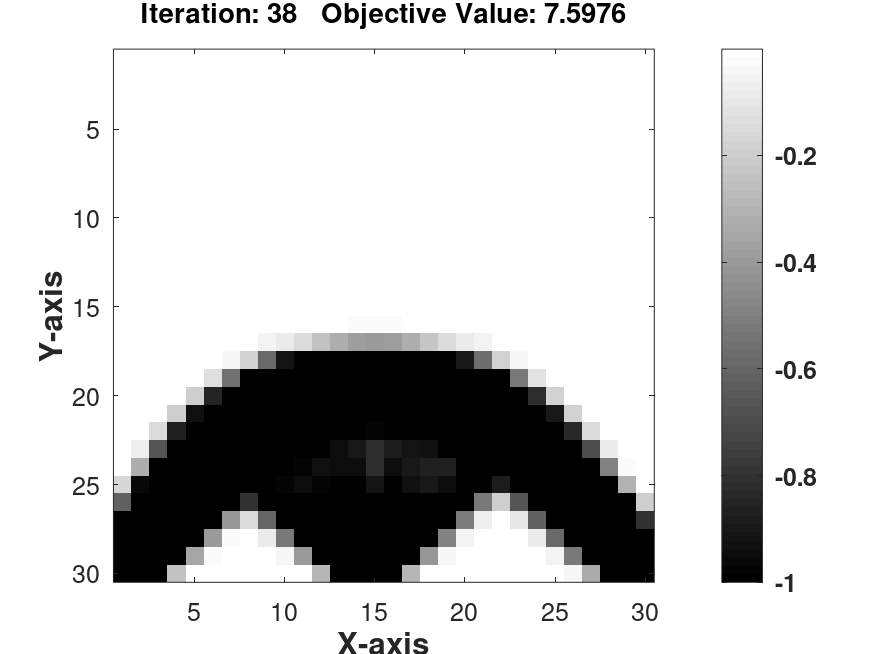} 
		Elements = 900  & 
		\includegraphics[width=5cm,height=4cm]{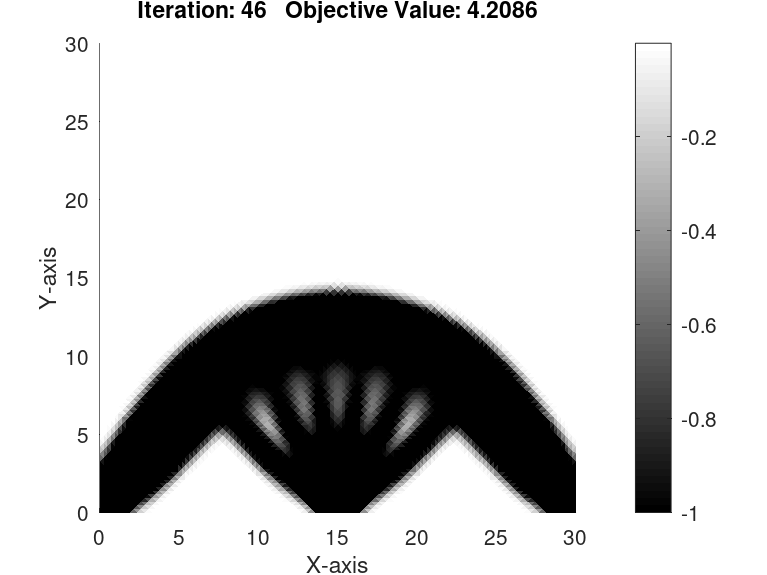} 
		Elements = 16384 &  
		\includegraphics[width=5cm,height=4cm]{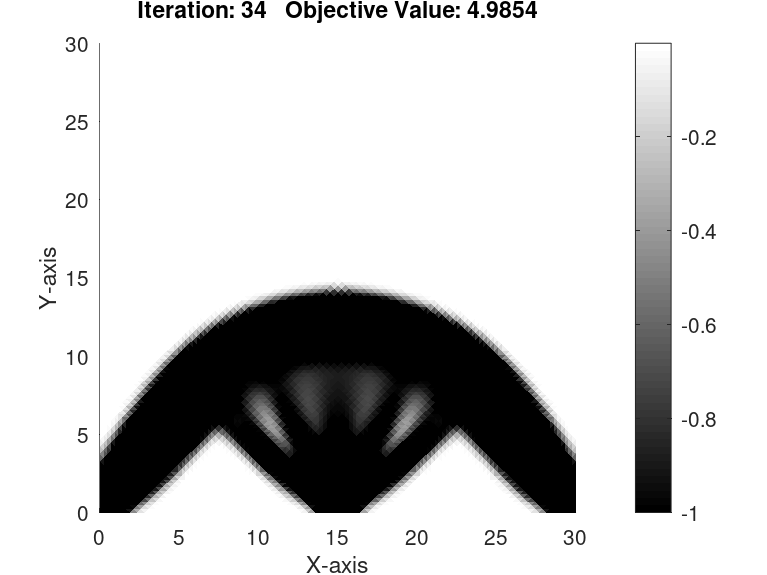}  
		Elements = 16384  \\
		\hline
		&
	 \includegraphics[width=5cm,height=4cm]{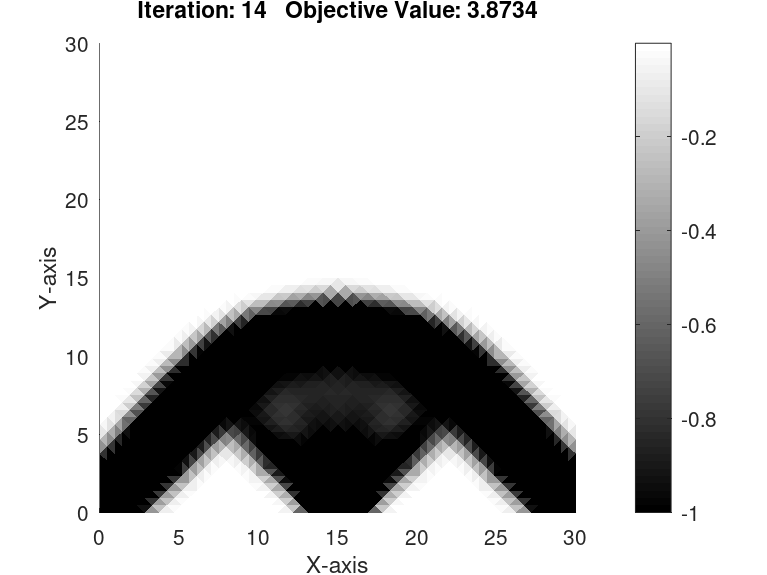} 
	 Elements = 4096 & 
	  \includegraphics[width=5cm,height=4cm]{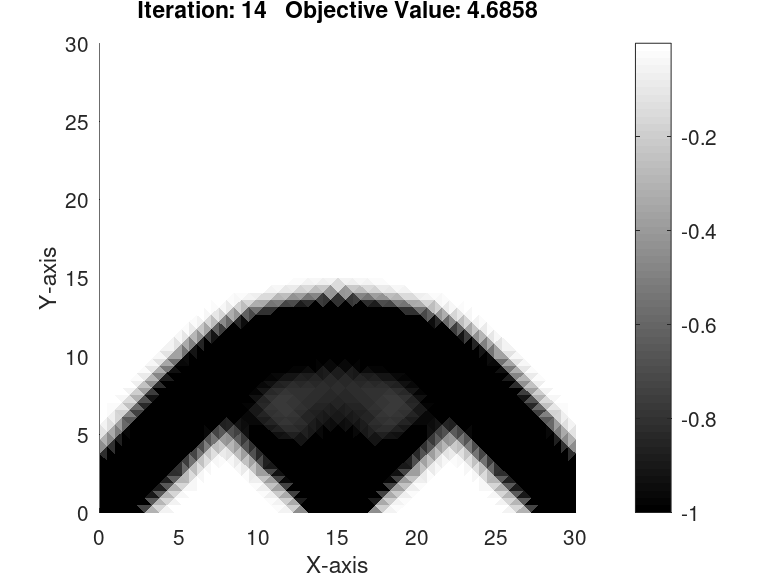} 
	 Elements = 4096
		 \\
		\hline

	\end{tabular}
	\caption{Comparison for \textbf{Bridge structure of Dimension $30 \times 30$}}
	\label{tab:bb1}
\end{figure}

	\begin{minipage}[t]{0.38\textwidth}  
		\centering
		\begin{tikzpicture}[baseline={(current bounding box.north)}, scale=1.0]
			\fill[gray!50] (0,0) rectangle (4,3);
			\draw[thick] (0,0) rectangle (4,3);
			\filldraw[black] (0,0) circle (0.05);
			\filldraw[black] (4,0) circle (0.05);
			\filldraw[red] (2,0) circle (0.07);
			\draw[very thick, ->, red] (2,-0.2) -- (2,-0.6);
			\node[below] at (2,-0.6) {\footnotesize Force};
			\filldraw[black] (0,0) -- (-0.2,-0.15) -- (0.2,-0.15) -- cycle;
			\filldraw[black] (4,0) -- (3.8,-0.15) -- (4.2,-0.15) -- cycle;
			\node[below] at (0,-0.3) {\footnotesize Fixed};
			\node[below] at (4,-0.3) {\footnotesize Fixed};
		\end{tikzpicture}
		
		
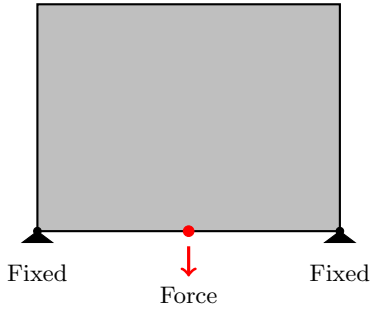
\captionof{figure}{\textbf{Bridge Structure} with fixed and force node. $E = 1$, $\nu = 0.3$, $volfrac = 0.3$, $penal = 3$, and $\vec{F} = (0, -1)$.}
		\label{fig:brbeam}
	\end{minipage}
	\hfill
	\begin{minipage}[t]{0.5\textwidth} 
	\textbf{Bridge Structure:} Figure~\ref{fig:brbeam} illustrates the details of the Bridge Structure, where the fixed nodes are located at the bottom two corners of the domain and a force is applied at the center of the bottom side, directed along the negative Y-axis. The parameters used are described within the same figure. Figure~\ref{tab:bb1} \& ~\ref{tab:bb2} presents the optimized topologies obtained using different resolutions.
\end{minipage}%

\begin{table}[htbp]
	\centering
	\renewcommand{\arraystretch}{1.2}
	
	\begin{tabular}{|
			>{\centering\arraybackslash}m{1.5cm}|
			>{\centering\arraybackslash}m{1.5cm}|
			>{\centering\arraybackslash}m{1.4cm}|
			>{\centering\arraybackslash}m{1.7cm}|
			>{\centering\arraybackslash}m{1.5cm}|
			>{\centering\arraybackslash}m{1.5cm}|
			>{\centering\arraybackslash}m{1.6cm}|
			>{\centering\arraybackslash}m{1.2cm}|
			>{\centering\arraybackslash}m{1.2cm}|}
			\hline
					
		\textbf{Element} &
		\textbf{No. of Elements} &
		\textbf{Final Objective Value} &
		\textbf{Iterations} &
		\textbf{Bulk Residual} &
		\textbf{Internal Jump Residual} &
		\textbf{Neumann Residual} &
		\textbf{Local Error $\eta^2$} &
		\textbf{Global Error $\eta$} \\
		
		\hline
		
		Q1 & 900 & 7.5976 & 38 & 4.7767 & 6.9280 & 0.6338 & 12.338 & 3.5125 \\
		\hline
		
		\multirow{2}{*}{P1}
		& 16384 & 4.2086 & 46 & 0 & 0.5003 & 0.5661 & 1.0663 & 1.0326 \\
		\cline{2-9}
		& 4096 & 3.8734 & 14 & 0 & 0.5669 & 0.7709 & 1.3378 & 1.1566 \\
		\hline
		
		\multirow{2}{*}{P2}
		& 16384 & 4.9854 & 34 & 0.9346 & 0.3814 & 0.0134 & 1.3296 & 1.1531 \\
		\cline{2-9}
		& 4096 & 4.6858 & 14 & 3.7298 & 0.3848 & 0.0942 & 4.2089 & 2.0516 \\
		\hline
		
	\end{tabular}
	
	\caption{Comparison of objective values, optimization iterations, and residual-based a posteriori error estimates for the bridge structure problem on a $30 \times 30$ domain using different finite element discretizations.}
	
	\label{tab:bridge_results}
	
\end{table}

	\begin{figure}[htbp]
	\centering
	\renewcommand{\arraystretch}{1.5}
	\begin{tabular}{|
			>{\centering\arraybackslash}m{5cm}|
			>{\centering\arraybackslash}m{5cm}|
			>{\centering\arraybackslash}m{5cm}|}
		\hline 
\textbf{$Q_1$} & \textbf{$P_1$}& \textbf{$P_2$} \\
\hline
		
		\includegraphics[width=5cm,height=5cm]{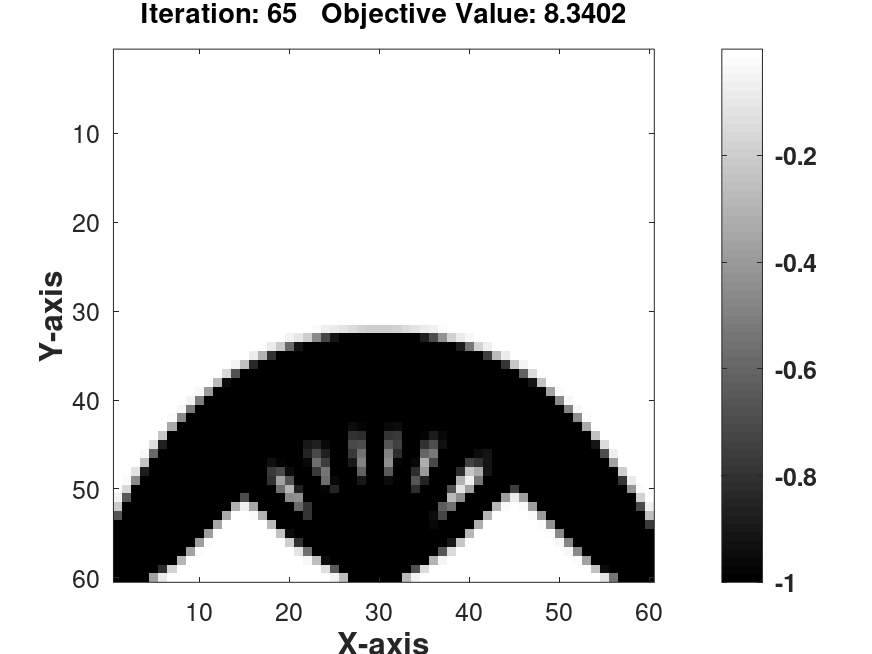} 
		Elements = 3600  & 
		\includegraphics[width=5cm,height=5cm]{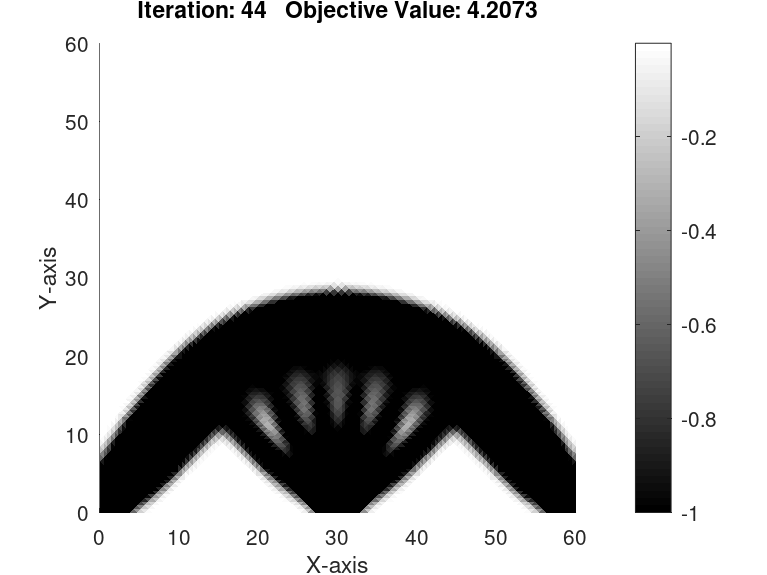} 
		Elements = 16384 &  
		\includegraphics[width=5cm,height=5cm]{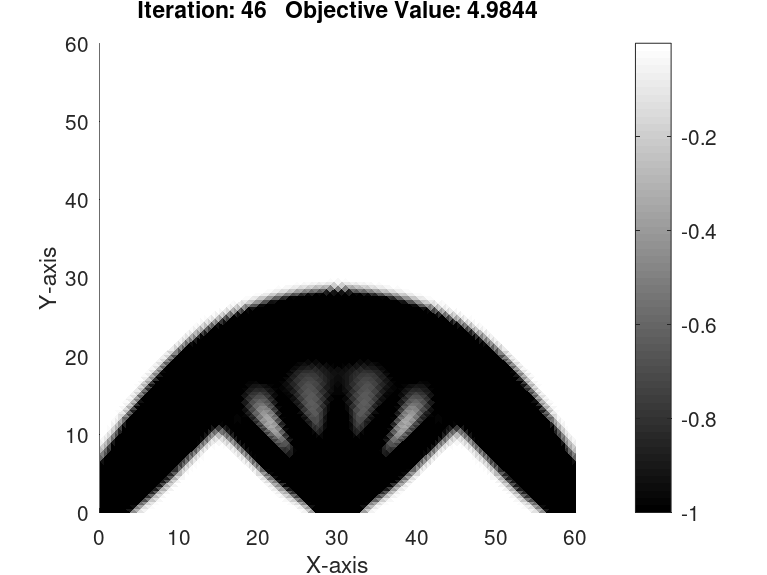}  
		Elements = 16384  \\
		\hline
		& 		\includegraphics[width=5cm,height=5cm]{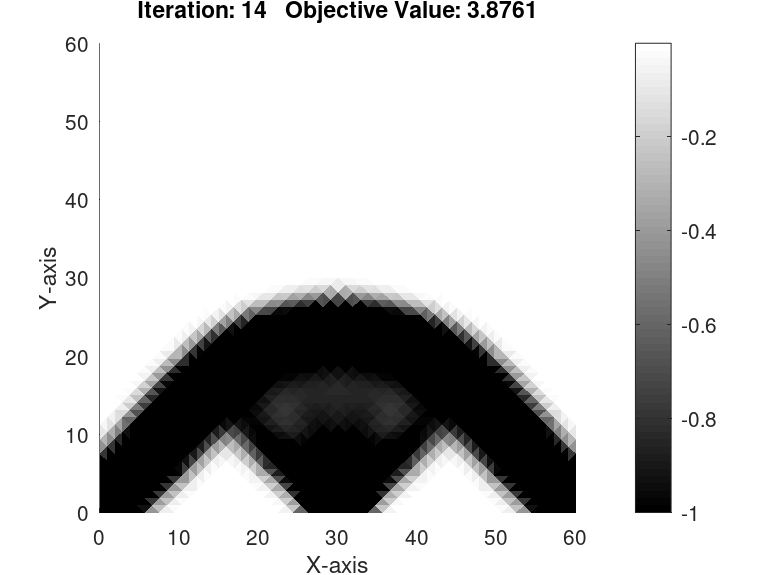} 
		Elements = 4096   
	 & 
	 \includegraphics[width=5cm,height=5cm]{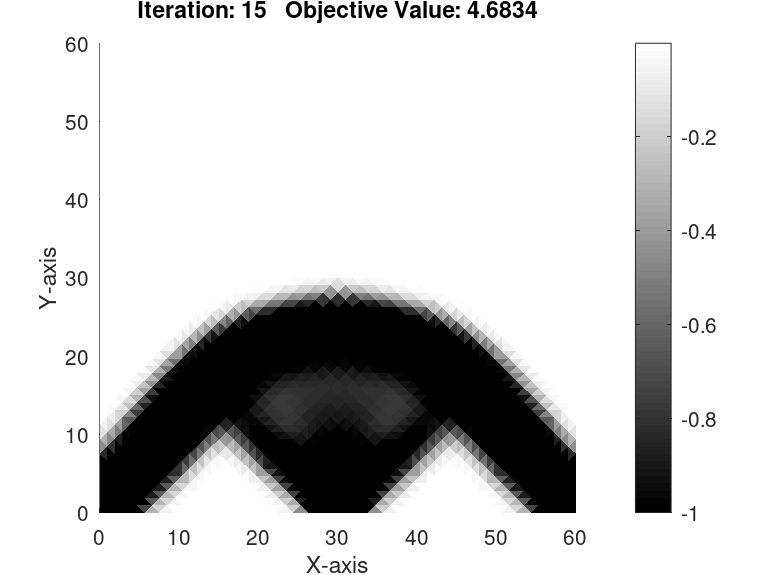} 
	 	Elements = 4096  
		\\
		\hline

	\end{tabular}
	\caption{Comparison for \textbf{Bridge structure of Dimension $60 \times 60$}}
	\label{tab:bb2}
\end{figure}

\begin{table}[H]
	\centering
	\renewcommand{\arraystretch}{1.2}
	
	\begin{tabular}{|
			>{\centering\arraybackslash}m{3cm}|
			>{\centering\arraybackslash}m{3cm}|
			>{\centering\arraybackslash}m{3cm}|
			>{\centering\arraybackslash}m{3cm}|}
		\hline 
		
		\textbf{Element Type} & 
		\textbf{No. of Elements} & 
		\textbf{Final Objective Value} & 
		\textbf{Iterations} \\
		\hline
		
		Q1 & 3600 & 8.3402 & 65 \\
		\hline
		
		\multirow{2}{*}{P1} &  16384 & 4.2073 & 44  \\
		\cline{2-4}
	     	& 4096 & 3.8761 & 14 \\
		\hline
		\multirow{2}{*}{P2} & 16384 & 4.9844 & 46 \\
		\cline{2-4}
	    	& 4096 & 4.6834 & 15 \\
		\hline
	\end{tabular}
	\caption{Bridge Structure of Domain dimension $60 \times 60$ }
	\label{tab:b1}
\end{table}

	\begin{minipage}[t]{0.45\textwidth}  

		\centering
		\begin{tikzpicture}[baseline={(current bounding box.north)}, scale=1]
			\draw[thick] (0,0) -- (4,1) -- (4,2) -- (0,3) -- cycle;
			\fill[gray!50] (0,0) -- (4,1) -- (4,2) -- (0,3) -- cycle;
			
			\foreach \y in {0.1,0.3,...,2.9}
			\draw[black] (-0.2,\y) -- (0,\y+0.1);
			
			\filldraw[red] (4,1.5) circle (0.05);
			\draw[very thick, ->, red] (4.2,1.5) -- (4.2,1.1);
			\node[right] at (4.2,1.1) {\footnotesize Force};
			
			\node[rotate=90] at (-0.5,1.5) {\footnotesize Fixed Nodes};
		\end{tikzpicture}
		
		\captionof{figure}{\textbf{Beveled Beam} with fixed and force node. $E = 1$, $\nu = 0.3$, $volfrac = 0.5$, $penal = 3$, and $\vec{F} = (0, -1)$.}
		\label{fig:bbbeam}
	\end{minipage}
	\hfill
		\begin{minipage}[t]{0.5\textwidth}  
		\textbf{Beveled Beam:} A practical beveled beam structure was considered for experimentation, where the beam is slanted, meaning that certain regions of the design domain are passive (i.e., inactive in the optimization). Figure~\ref{fig:bbbeam} illustrates the configuration of the beveled beam, where the fixed nodes are located entirely along the west side of the domain. A downward force is applied at the mid-point node on the east boundary, acting along the negative Y-axis. 
	\end{minipage}%

	\begin{figure}[H]
	\centering
	\renewcommand{\arraystretch}{1.5}
	\begin{tabular}{|
			>{\centering\arraybackslash}m{5cm}|
			>{\centering\arraybackslash}m{5cm}|
			>{\centering\arraybackslash}m{5cm}|}
		\hline 
\textbf{$Q_1$} & \textbf{$P_1$}& \textbf{$P_2$} \\
\hline
		
		\includegraphics[width=5cm,height=4cm]{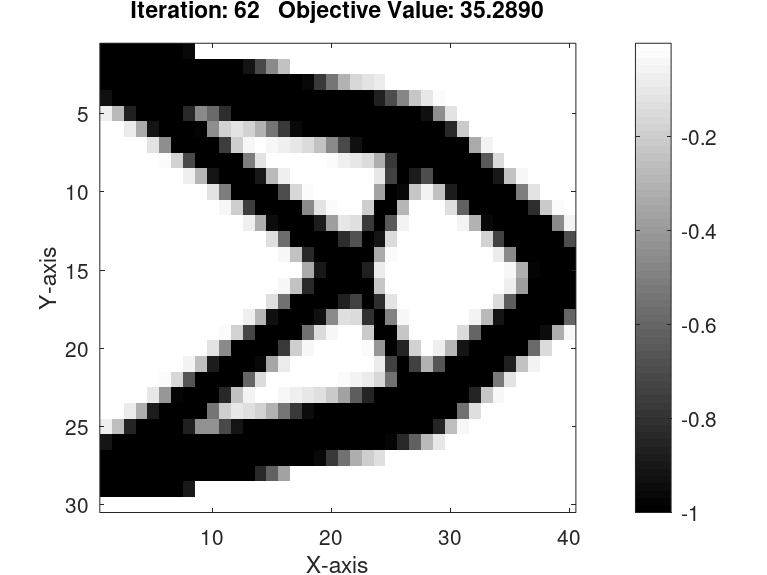} 
		Elements = 1200  & 
		\includegraphics[width=5cm,height=4cm]{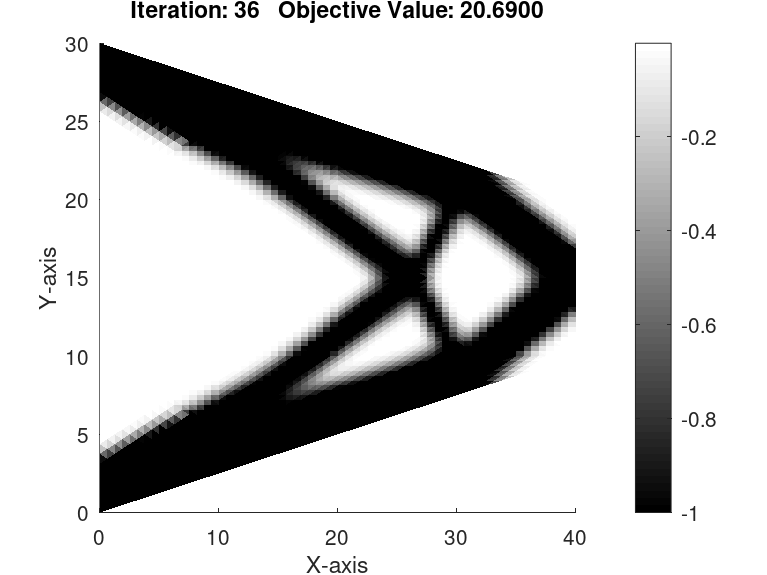} 
		Elements = 7424 &  
		\includegraphics[width=5cm,height=4cm]{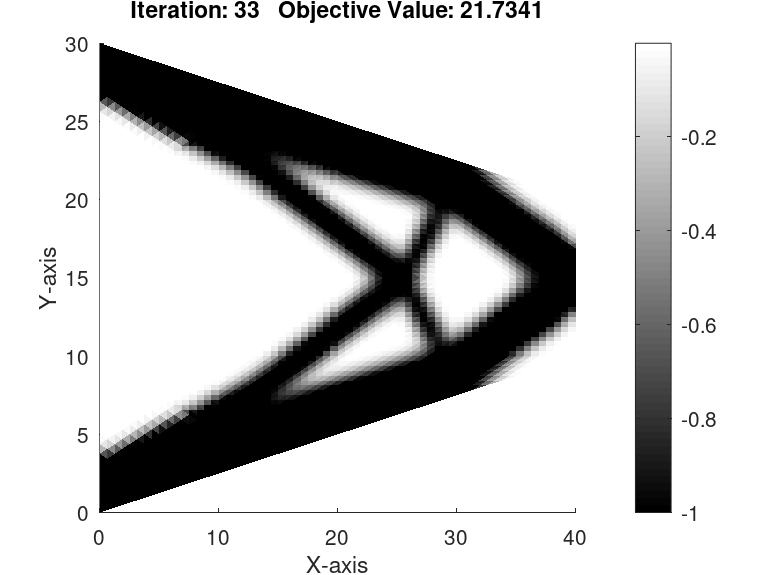}  
		Elements = 7424  \\
		\hline

	\end{tabular}
	\caption{Comparison for \textbf{Beveled Beam of Dimension $40 \times 30$}}
	\label{tab:bbeam2}
\end{figure}

The optimized layouts obtained across the cantilever beam, bridge structure, and beveled beam problems exhibit noticeable variations depending on the finite element discretization used. In general, the topology generated using $Q_1$ bilinear quadrilateral elements tends to produce relatively thicker and simpler load-carrying members, reflecting the limited geometric flexibility of quadrilateral meshes. In contrast, the solutions obtained using $P_1$ linear triangular elements form more intricate structural networks with finer members and additional branching patterns, enabling a better representation of complex stress transfer mechanisms within the design domain. The use of $P_2$ quadratic triangular elements further improves the smoothness and continuity of the structural boundaries, resulting in more refined load paths and a more realistic representation of stress distributions. These observations are consistent across all considered benchmark problems and demonstrate that triangular discretizations provide greater geometric flexibility and improved capability to capture detailed structural features, thereby leading to a richer variety of optimized topologies compared to $Q_1$ elements.

For all benchmark problems, the $Q_1$ discretization produces the highest compliance values, indicating comparatively less stiff optimized designs. In addition, the global error estimator $\eta$ associated with $Q_1$ elements is significantly larger than those obtained using triangular discretizations. While assessing the accuracy the examination of the a-posteriori error components shows that the internal jump residuals contribute most strongly to this estimate, suggesting reduced continuity of normal stresses across element
interfaces. This indicates that bilinear quadrilateral elements may struggle to accurately capture stress gradients in regions where the optimized topology develops thin members or complex load paths. The $P_1$ triangular discretization consistently achieves lower compliance values across all three structures, indicating improved structural stiffness. Moreover, the corresponding global error estimates are substantially smaller than those of $Q_1$ elements. The bulk residuals for $P_1$ meshes are zero due to second derivation of linear elements. The dominant contributions arise from internal jump and Neumann boundary residuals, which remain moderate compared to the quadrilateral case. These observations confirm that linear triangular elements provide a more accurate representation of the displacement field while maintaining computational efficiency. Quadratic $P_2$ elements further improve the smoothness of the displacement and stress fields due to their higher-order interpolation capability. This generally leads to reduced jump residuals on sufficiently refined meshes. However, the bulk residual component may increase on coarser meshes, reflecting the sensitivity of higher-order approximations to mesh resolution. Although $P_2$ elements produce physically smoother topologies, they do not consistently yield lower compliance than $P_1$ elements and require substantially greater computational effort. The bridge and beveled beam problems exhibit trends similar to those observed for the cantilever beam. In all cases, triangular
discretizations outperform quadrilateral elements in both compliance and error estimation. The consistency of these results across different geometries, boundary conditions, and loading scenarios demonstrates the robustness of triangular elements for topology optimization.
	
\section{Conclusion}

This study investigated the influence of finite element discretization on topology optimization results for three benchmark problems: the cantilever beam, bridge structure, and beveled beam. The optimized layouts obtained from different discretizations exhibit noticeable structural variations. The $Q_1$ bilinear quadrilateral elements generally produce thicker and simpler load-carrying members, whereas triangular discretizations allow the formation of finer and more intricate structural networks. The $P_2$ quadratic triangular elements further improve the smoothness of structural boundaries due to their higher-order interpolation capability. The numerical results indicate that the $Q_1$ discretization consistently produces the largest global a posteriori error estimator $\eta$ and higher compliance values across all benchmark problems. In contrast, triangular discretizations provide improved accuracy and structural performance. In particular, the $P_1$ triangular elements achieve the lowest compliance values while maintaining relatively small error estimates, demonstrating an effective balance between accuracy and computational efficiency. Although $P_2$ elements yield smoother displacement and stress fields, their computational cost is higher and they do not consistently outperform $P_1$ elements in terms of compliance.

Overall, the results demonstrate that triangular finite element discretizations, particularly $P_1$, offer a more reliable and accurate framework for SIMP-based topology optimization when assessed using residual-based a posteriori error estimation. These findings highlight the importance of element choice in achieving accurate stress representation and efficient structural designs in topology optimization problems.

\section*{Conflicts of Interest}

On behalf of all authors, the corresponding author declares that there are no conflicts of interest related to the publication of this work. 

\section*{Acknowledgments}
The authors would also like to acknowledge Dr. Tapan Kumar Hota for his insightful remarks and valuable observations on this paper.

\section*{Data Availability Statement}
The data presented in this study are available on request from the corresponding author due to some of the data involves privacy.	

\section*{Grammar and Readability Disclosure }
This document has been reviewed with AI-based tools \cite{chatgpt2023} to check grammar and readability improvements. 
	
	\bibliographystyle{unsrt}
	\bibliography{reference}

\end{document}